\documentclass[11pt,a4paper]{article}
\usepackage{latexsym,amssymb,amsthm,amsmath}
\title{Optimal Hardy-Sobolev-Maz'ya inequalities with multiple
 interior singularities}

\author{\Large Stathis Filippas$^{1,4}$,  Achilles  Tertikas$^{2,4}$
\& Jesper Tidblom$^{3}$  \\
                                                                           \\
        Department of Applied Mathematics$^{1}$ \\
         University of Crete,
         71409 Heraklion,  Greece \\
        filippas@tem.uoc.gr\\
                  \\
  Department of Mathematics$^{2}$ \\
     University of Crete,
         71409 Heraklion,  Greece          \\
  tertikas@math.uoc.gr\\
\\
    The Erwin Schr\"odinger Institute (ESI)$^{3}$ \\
          Boltzmanngasse 9, A-1090 Vienna, Austria \\
          Jesper.Tidblom@esi.ac.at\\
                                      \\
        Institute of Applied and Computational Mathematics$^4$, \\
        FORTH, 71110 Heraklion, Greece \\
    \\ }

\date{}
\newtheorem{theorem}{Theorem}[section]
\newtheorem{proof*}{Proof:}
\newtheorem{corollary}[theorem]{Corollary}
\newtheorem{lemma}[theorem]{Lemma}
\numberwithin{equation}{section}
\newcommand{\be}{\begin{equation}}
\newcommand{\ee}{\end{equation}}
\newcommand{\bea}{\begin{eqnarray}}
\newcommand{\eea}{\end{eqnarray}}
\newcommand{\la}{\label}
\newcommand{\xa}{\alpha}
\newcommand{\xb}{\beta}
\newcommand{\xg}{\gamma}
\newcommand{\xe}{\varepsilon}
\newcommand{\xs}{\sigma}

\newcommand{\ra}{\rightarrow}
\newcommand{\rn}{\mathbb{R}^n}

\def\finedim{$\hfill \Box$}

\begin{document}
\maketitle

 \begin{center} {\bf  Dedicated to Prof. Vladimir  Maz'ya with esteem}
\end{center}
 \vspace{.5cm}

\begin{abstract}
 \noindent
In this article
 we first establish  a complete characterization
of Hardy's inequalities in $\mathbb{R}^n$
 involving distances
to different codimension subspaces.  In particular the corresponding potentials
 have strong interior singularities.
 We then provide  necessary and sufficient conditions
for the validity of  Hardy-Sobolev-Maz'ya inequalities  with optimal  Sobolev terms.
\end{abstract}


\section{Introduction}
For $n \geq 3$ we write $\mathbb{R}^n =\mathbb{R}^k \times \mathbb{R}^{n-k}$, with $1 \leq k \leq n$. We also
introduce the codimension $k$   affine subspace
\[
S_k:=\{x=(x_1, \dots  x_k, \ldots x_n) \in \mathbb{R}^n:~ x_1=\ldots=x_k=0\}.
\]
The Euclidean  distance of a point $x \in  \mathbb{R}^n$ from $S_k$ is then  given by
\[
d(x)=d(x, S_k)=  |  \mathbf{X_k} |,~~~~ {\rm where}~~~~\mathbf{X_k} := (x_1,\ldots,x_k,0,\ldots,0).
\]
The classical Hardy inequality in $\mathbb{R}^n$ when distance is taken from $S_k$, reads
\be\la{0}
 \int_{\mathbb{R}^n}|\nabla u|^2dx \geq  \left( \frac{k-2}{2} \right)^2
 \int_{\rn}\frac{u^2}{|\mathbf{X_k} |^2}dx, ~~~~~~~~~~~ u \in
 C^{\infty}_0(\rn \setminus S_k),
\ee
where the constant
$\frac{(k-2)^2}{4}$ is the optimal one. This result has been
improved and generalized in many different ways, see for example \cite{ACR, AE, An, BFT, BM, BV, CF, Dav,
 DELV,  FT, Hoff2,  Tid1, Tid2} and references therein.

On the other hand the standard Sobolev inequality with critical exponent  states that
\[
   \int_{\rn}|\nabla u|^2dx  \geq S_n
 \left( \int_{\rn}|u|^{\frac{2n}{n-2}}dx \right)^{\frac{n-2}{n}},
   \quad u \in C_0^{\infty}(\rn),
\]
 where  $S_n=
\pi n(n-2) \left( \frac{\Gamma(\frac{n}{2})}{\Gamma(n)}
 \right)^{2/n}$ is   the best  Sobolev  constant, see \cite{A, T}. For versions of Sobolev
inequalities involving subcritical exponents and weights see e.g.
\cite{AFT2,  BT, CKN}.

Maz'ya, in his book, combined both inequalities when $1\leq k \leq n-1$, establishing
that for any $ u \in
 C^{\infty}_0(\rn \setminus S_k)$
\be\label{maz}
 \int_{\mathbb{R}^n}|\nabla u|^2dx  \geq  \left( \frac{k-2}{2} \right)^2
 \int_{\rn}\frac{u^2}{|\mathbf{X_k} |^2}dx +  c_{k,Q}
 \left(\int_{\mathbb{R}^n}|\mathbf{X_k} |^{\frac{Q-2}{2}n-Q}~|u|^{Q}dx
  \right)^{\frac{2}{Q}},
\ee
for $2<Q\leq 2^{*}=\frac{2n}{n-2}$; cf.  \cite{Maz}, Section 2.1.6/3. Concerning the best constant
$c_{k,Q}$, it was shown  in \cite{TT} that $c_{k,2^*}< S_n$  for $3 \leq k \leq n-1$, $n \geq   4$  or $k=1$ and $n \geq 4$.
Surprisingly, in the case   $k=1$ and  $n=3$ Benguria Frank and Loss  \cite{BFL}  (see also
Mancini and Sandeep \cite{MS}) established that $c_{1,6}= S_3=3 (\pi/2)^{4/3}$!
Maz'ya and Shaposhnikova \cite{Maz2}  have recently computed the best constant in the  case $k=1$ and $Q=\frac{2(n+1)}{n-1}$.
These are the only cases where the best constant $  c_{k,Q} $ is known.  For other type of Hardy--Sobolev
inequalities see \cite{FMT1, Fil, PT}.

In case  $k=n$, that is,  when distance is taken from the origin, inequality (\ref{maz}) fails.  Brezis and Vazquez
\cite{BV} considered a bounded domain containing the origin and improved the Hardy inequality by adding
a subcritical Sobolev term.   It turns out that in a bounded domain one can have the critical Sobolev exponent
at the expense however of adding a logarithmic weight. More specifically let
\[
X(t) = (1- \ln t)^{-1}, ~~~~~~~~~0< t <1.
\]
Then the analogue of (\ref{maz}) in the case of a bounded domain $\Omega$
 containing the origin,
for the  critical exponent reads:
\be\la{thb}
\int_{\Omega} |\nabla u|^2 dx -   \left(\frac{n-2}{2} \right)^2   \int_{\Omega} \frac{u^2}{|x|^2} dx
 \geq
C_n(\Omega) \left( \int_{\Omega} X^{\frac{2(n-1)}{n-2}} \left(\frac{|x|}{D} \right) |u|^{\frac{2n}{n-2}}  dx
 \right)^{\frac{n-2}{n}},~~~u \in C^{\infty}_0(\Omega),
\ee
where  $D=\sup_{x \in \Omega}|x|$; cf \cite{FT}. The best constant in (\ref{thb}) was recently computed in
\cite{AFT} and is given by
\[
C_n(\Omega) = (n-2)^{-\frac{2(n-1)}{n}} S_n.
\]
It is worth noticing that in the case $n=3$ once again one has $C_3(\Omega)=S_3=3 (\pi/2)^{4/3}$!

In a recent work  \cite{Fil3} we studied Hardy--Sobolev--Maz'ya inequalities  that involve distances  taken
from different codimension subspaces of the boundary.  In particular, working in the upper half space
$ \mathbb{R}^n_+=\{ x \in \mathbb{R}^n: x_1 >0 \} $
  and taking distances  from $S_k \subset \partial  \mathbb{R}^n_+ \equiv  S_1$,  $k=1,2, \ldots,n$,
we have established  that the following inequality holds true for any $u \in C_0^{\infty}(\mathbb{R}^n_+)$
\be\label{0.3}
  \int_{\mathbb{R}^n_+}|\nabla u|^2dx \geq \int_{\mathbb{R}^n_+}\left(\frac{\beta_1}{x_1^2}+
 \frac{\beta_2}{|\mathbf{X_2} |^2} + \ldots + \frac{\beta_n}{|\mathbf{X_n} |^2} \right)u^2dx,
\ee
if and only if there exist  nonpositive constants
 $\xa_1, \ldots, \xa_n$, such that
\begin{eqnarray}\la{bs}
 \beta_1 = -\alpha_1^2   +\frac14, ~~~~~
 \beta_m  =  - \alpha_m^2 + \left(\alpha_{m-1}-\frac12 \right)^2,
~~~~m=2,3,\ldots,n.
\end{eqnarray}
 Moreover if $\alpha_n <0$ one can add in the right hand side the critical Sobolev term,
thus obtaining  the  Hardy--Sobolev--Maz'ya inequality valid for any $u \in C_0^{\infty}(\mathbb{R}^n_+)$
\be\label{0.4}
  \int_{\mathbb{R}^n_+}|\nabla u|^2dx \geq \int_{\mathbb{R}^n_+}\left(\frac{\beta_1}{x_1^2}+
 \frac{\beta_2}{|\mathbf{X_2} |^2} + \ldots + \frac{\beta_n}{|\mathbf{X_n} |^2} \right)u^2dx
+  C\left(\int_{\mathbb{R}^n_+}|u|^{\frac{2n}{n-2}}dx
  \right)^{\frac{n-2}{n}};
\ee
we refer to  \cite{Fil3} for the detailed statements.

In the present work we  consider the case where  distances are again  taken from  different codimension subspaces
 $S_k \subset \mathbb{R}^n$,
which however are now  placed in the interior of the domain  $\mathbb{R}^n$. We consider the cases  $k=3,\ldots, n$
 since there is no positive Hardy constant in case $k=2$ (cf (\ref{0}))  and the case $k=1$ corresponds to the case studied in
 \cite{Fil3}.

More precisely our first result reads

\quad
\\ \textbf{Theorem A}
     ({\bf Improved Hardy inequality}) \\
 {\it Suppose  $n \geq 3$. \\
  {\bf i)}Let   $\alpha_3, \alpha_4,
\ldots, \alpha_n$  be arbitrary real numbers and
\begin{eqnarray*}
 \beta_3 = -\alpha_3^2   +\frac14,  ~~~~~
 \beta_m  =  - \alpha_m^2 + \left(\alpha_{m-1}-\frac12 \right)^2,
~~~~m=4,\ldots,n.
\end{eqnarray*}
Then for any $u \in C_0^{\infty}(\mathbb{R}^n)$ there holds
\[
  \int_{\mathbb{R}^n}|\nabla u|^2dx \geq \int_{\mathbb{R}^n}\left(
 \frac{\beta_3}{|\mathbf{X_3} |^2} + \ldots + \frac{\beta_n}{|\mathbf{X_n}|^2} \right)u^2dx.
\]
{\bf ii)}Suppose that for some real numbers $\beta_3, \beta_4
\ldots, \beta_n$ the
 following
inequality holds \be\nonumber
  \int_{\mathbb{R}^n}|\nabla u|^2dx \geq
\int_{\mathbb{R}^n}\left(
 \frac{\beta_3}{|\mathbf{X_3} |^2} + \ldots
 + \frac{\beta_n}{|\mathbf{X_n} |^2} \right)u^2dx,
\ee for any $u \in C_0^{\infty}(\mathbb{R}^n)$. Then, there
exists  nonpositive constants
 $\xa_3, \ldots, \xa_n$, such that
\begin{eqnarray*}
 \beta_3 = -\alpha_3^2   +\frac14,  ~~~~~
 \beta_m  =  - \alpha_m^2 + \left(\alpha_{m-1}-\frac12 \right)^2,
~~~~m=4,\ldots,n.
\end{eqnarray*}
}

We note that the recursive formula for the $\beta$'s  in the above Theorem,   is the same as in (\ref{bs}).
However,  since the coefficients in the above Theorem start from $ \beta_3$ -- and not  from  $ \beta_1$--
the  best  constants in the case of interior singularities  are different from
the  best  constants  when singularities  of the same codimension are placed on the boundary.
See for instance Corollary \ref{co2.3} and compare with  Corollary  2.4  of \cite{Fil3}.

To state our next results we define
\begin{eqnarray}\la{betas}
 \beta_3 = -\alpha_3^2   +\frac14,   ~~~~~~~
 \beta_m  =  - \alpha_m^2 + \left(\alpha_{m-1}-\frac12 \right)^2,
~~~~m=4,\ldots,n.
\end{eqnarray}
Our next theorem gives  a complete answer  as to when  we can add a Sobolev term.

\quad \\ \textbf{Theorem B}
({\bf  Improved Hardy--Sobolev--Maz'ya  inequality}) \\
 {\it Let $\alpha_3, \alpha_4, \ldots,
\alpha_n$, $n \geq 3$,  be arbitrary nonpositive
 real numbers and  $ \beta_3, \ldots, \beta_n$ are given  by (\ref{betas}).
Then,  if $\xa_n <0$ there exists a positive constant $C$ such
that
 for any
$u \in C_0^{\infty}(\mathbb{R}^n)$ there holds
\be
\la{1.20}
  \int_{\mathbb{R}^n}|\nabla u|^2dx \geq \int_{\mathbb{R}^n}\left(
 \frac{\beta_3}{|\mathbf{X_3} |^2} + \ldots + \frac{\beta_n}{|\mathbf{X_n}|^2} \right)u^2dx
 + C\left(\int_{\mathbb{R}^n} |\mathbf{X_2} |^{\frac{Q-2}{2}n-Q}|u|^{Q}dx
  \right)^{\frac{2}{Q}},
\ee
for any $2 < Q \leq \frac{2n}{n-2}$.

If  $\xa_n=0$ then there is no positive constant $C$ such that
(\ref{1.20}) holds. }
\quad \\

The  above  result extends considerably the original inequality  by Maz'ya (\ref{maz}). First
by having at the same time,
  all possible combinations  of Hardy potentials involving the  distances $|\mathbf{X_3} |, \ldots ,|\mathbf{X_n} |$.
In addition the weight in the Sobolev term is  stronger than the weight used in (\ref{maz}).

We note that a similar result can be produced in the setting of \cite{Fil3} where singularities are placed
on the boundary $\partial \mathbb{R}^n_+$. More precisely the following  inequality holds true for
any $u \in C_0^{\infty}(\mathbb{R}^n_+)$
\be
\la{1.20b}
  \int_{\mathbb{R}^n_+}|\nabla u|^2dx \geq \int_{\mathbb{R}^n_+}\left(
 \frac{\beta_1}{x_1^2} +  \frac{\beta_2}{|\mathbf{X_2}|^2}  \ldots + \frac{\beta_n}{|\mathbf{X_n}|^2} \right)u^2dx
 + C\left(\int_{\mathbb{R}^n_+} x_1^{\frac{Q-2}{2}n-Q}|u|^{Q}dx
  \right)^{\frac{2}{Q}},
\ee
provided that $\alpha_n <0$,
where the constants $\beta_i$ are given by (\ref{bs}) and  $2 < Q \leq \frac{2n}{n-2}$.
In this case the weight in the right hand side is even stronger than the one in (\ref{1.20}).
In the light of  (\ref{1.20b}) one may ask whether one can replace the weight $|\mathbf{X_2}|$ in (\ref{1.20})
by $|x_1|$.  It turns out that this is possible provided we properly restrict the exponent $Q$. More
precisely we have:

\quad \\ \textbf{Theorem C}
({\bf Improved Hardy--Sobolev--Maz'ya  inequality}) \\
 {\it Let $\alpha_3, \alpha_4, \ldots,
\alpha_n$, $n \geq 3$,  be arbitrary nonpositive
 real numbers  and  $ \beta_3, \ldots, \beta_n$ are given  by (\ref{betas}).
Then,  if $\xa_n <0$ there exists a positive constant $C$ such
that
 for any
$u \in C_0^{\infty}(\mathbb{R}^n)$ there holds
\be
\la{1.20a}
  \int_{\mathbb{R}^n}|\nabla u|^2dx \geq \int_{\mathbb{R}^n}\left(
 \frac{\beta_3}{|\mathbf{X_3} |^2} + \ldots + \frac{\beta_n}{|\mathbf{X_n}|^2} \right)u^2dx
 + C\left(\int_{\mathbb{R}^n} |x_1|^{\frac{Q-2}{2}n-Q}|u|^{Q}dx
  \right)^{\frac{2}{Q}},
\ee
for any $\frac{2(n-1)}{n-2} < Q \leq \frac{2n}{n-2}$.

If  $\xa_n=0$ then there is no positive constant $C$ such that
(\ref{1.20a}) holds. }

\quad \\

It is easily seen that the range of the exponent $Q$ in Theorem C is optimal since otherwise
the weight is not locally integrable.  In the special case $\beta_3=\ldots \beta_n=0$, the
corresponding weighted Sobolev inequality in (\ref{1.20a}) was proved by Maz'ya, cf \cite{Maz}
section 2.1.6/2.

An important role in our analysis  is played by two weighted
Sobolev inequalities, which are of independent interest; see
Theorems \ref{3.1} and \ref{3.2}.

The paper is organized as follows. In section 2 we give the proof
of Theorem A. In section 3 we give the proofs of Theorems B and  C.
 The  main ideas   are similar to the ones used in
\cite{Fil3} to  which we refer on various occasions.  On the other
hand ideas or technical estimates that are different from
\cite{Fil3}  are presented in detail.

\medskip
\noindent {\bf Acknowledgments} JT is thanking the
Departments of Mathematics and Applied Mathematics of University
of Crete for the invitation as well as the warm hospitality.

\section{Improved Hardy inequalities with multiple singularities}

The following simple lemma may be found in \cite{Fil3}.

\begin{lemma} \label{mainthm}
(i) Let $\mathbf{F} \in C^1(\Omega)$, then
\begin{equation} \label{eq1}
\int_{\Omega}|\nabla u|^2 dx = \int_{\Omega}\left( \mathrm{div}
\mathbf{F} - |\mathbf{F}|^{2}\right)|u|^2dx+\int_{\Omega} |\nabla
u+\mathbf{F}u|^2dx, ~~~~ \forall u \in C^{\infty}_0(\Omega).
\end{equation}
(ii) Let $\phi >0$, $\phi \in C^2(\Omega)$ and $u=\phi v$, then we
have
\begin{equation} \label{eq2}
\int_{\Omega}|\nabla u|^2dx = -\int_{\Omega}\frac{\Delta
\phi}{\phi}u^2dx + \int_{\Omega} \phi^2 |\nabla v|^2dx, ~~~~
\forall u \in C^{\infty}_0(\Omega).
\end{equation}
\end{lemma}
\begin{proof}
By expanding the square we have
\[
\int_{\Omega} |\nabla u+\mathbf{F}u|^2dx =\int_{\Omega}|\nabla
u|^2 dx +\int_{\Omega}|\mathbf{F}|^2u^2dx+\int_{\Omega}\mathbf{F}
\cdot \nabla u^2dx.
\]
Identity (\ref{eq1}) now follows by integrating by parts the
last term.

To prove (\ref{eq2}) we apply (\ref{eq1}) to
$\mathbf{F}=-\frac{\nabla \phi}{\phi}$. Elementary calculations
now yield the result.

\end{proof}

Let us recall our notation
\[
\mathbf{X_k} := (x_1,\ldots,x_k,0,\ldots,0) \quad \textrm{so~that }
\quad   |  \mathbf{X_k} |^2  =x_1^2 + \ldots + x_k^2;
\]
in particular $ |  \mathbf{X_n} |=|x|$.
We now give the proof of the first part of Theorem A:\\

\noindent
{\em Proof of Theorem A part (i):}
Let  $\gamma_3$, $\gamma_4$, $\ldots$, $\gamma_n$ be arbitrary real numbers
and put
\[
 \phi:=|\mathbf{X_3}|^{-\gamma_3}|\mathbf{X_4}|^{-\gamma_2}\cdot \ldots \cdot
 |\mathbf{X_n}|^{-\gamma_n},
\]
and
\[
 \mathbf{F}: = - \frac{\nabla \phi}{\phi}.
\]
An easy  calculation shows that
\[
 \mathbf{F} = \sum_{m=3}^n  \gamma_m   \frac{\mathbf{X_m}  }
{ |\mathbf{X_m}|^2}.
\]
With this choice of $\mathbf{F}$, we get
\[
\mathrm{div} \mathbf{F} = \sum_{m=3}^n   \gamma_m
  \frac{(m-2)}{ |\mathbf{X_m}|^2 },
\]
and
\[
|\mathbf{F} |^2 = \sum_{m=3}^n \frac{\gamma_m^2}{ |\mathbf{X_m}|^2 }
  + 2 \sum_{m=3}^n  \sum_{j=1}^{m-1}
\gamma_m \gamma_j  \frac{\mathbf{X_m}  }
{ |\mathbf{X_m}|^2}  \frac{\mathbf{X_j}  }
{ |\mathbf{X_j}|^2}
=
 \sum_{m=3}^n
\frac{\gamma_m^2}{ |\mathbf{X_m}|^2 }
+
 2 \sum_{m=3}^n  \sum_{j=1}^{m-1}
 \frac{\gamma_m \gamma_j}
 {|\mathbf{X_j}|^2  }.
\]
We then get that
\be\la{2.5}
-\frac{\Delta \phi}{\phi} =
\mathrm{div} \mathbf{F} - |\mathbf{F}|^2 = \sum_{m=3}^n
\frac{\beta_m}{|\mathbf{X_m}|^2},
\ee
where
\begin{eqnarray*}
 \beta_3 &=& -\gamma_3(\gamma_3-1),  \\
 \beta_m  &=&  - \gamma_m( 2-m +  \gamma_m  + 2 \sum_{j=3}^{m-1} \gamma_j ),
~~~~m=4,5,\ldots,n.
\end{eqnarray*}
We  next  set
\begin{eqnarray*}
 \gamma_3 &=& \alpha_3 + \frac12,  \\
 \gamma_m  &=& \alpha_m- \alpha_{m-1} + \frac12,
~~~~m=4,5,\ldots,n.
\end{eqnarray*}
With this choice of $\gamma$'s
 the $\beta$'s are given as in the statement of the Theorem.

We will use  Lemma \ref{mainthm}  with $\Omega =\mathbb{R}^n \setminus K_3$, where $K_3:= \{x \in    \mathbb{R}^n:
x_1=x_2=x_3=0 \}$. We have
\be\la{2.6}
 \int_{\mathbb{R}^n}|\nabla u|^2 dx \geq
 \int_{\mathbb{R}^n}\left(\mathrm{div}\mathbf{F}
  - |\mathbf{F} |^2\right) u^2  dx, ~~~~~~~~u \in C_0^{\infty}(\mathbb{R}^n \setminus K_3).
\ee
By a standard density argument  (\ref{2.6})  is  true even for $u \in C_0^{\infty}(\mathbb{R}^n)$.
The result then  follows from (\ref{2.5}) and (\ref{2.6}).

\finedim

\noindent


Some  interesting  cases   are presented  in  the following corollary.
\begin{corollary} \label{mainineq}
Let k=3,\ldots,n, $n \geq 3$,     and $ u \in C_0^{\infty}(\mathbb{R}^n)$.  Then
\begin{eqnarray}\la{cor1}
   \int_{\mathbb{R}^n}|\nabla u|^2dx \geq \int_{\mathbb{R}^n}\left(
  \left( \frac{k-2}{2} \right)^2 \frac{1}{|\mathbf{X_k}|^2} + \frac{1}{4}\frac{1}{|\mathbf{X_{k+1}}|^2}
 \ldots  \right.
  + \left. \frac{1}{4}\frac{1}{|\mathbf{X_n}|^2} \right)u^2dx.
\end{eqnarray}
Also,
\begin{eqnarray}\la{cor2}
   \int_{\mathbb{R}^n}|\nabla u|^2 dx \geq
  \left( \frac{k-2}{2} \right)^2 \int_{\mathbb{R}^n}   \frac{u^2}{|\mathbf{X_k}|^2}  dx
+   \left( \frac{n-k}{2} \right)^2 \int_{\mathbb{R}^n}  \frac{u^2}{|x|^2} dx.
\end{eqnarray}
\end{corollary}
\begin{proof}We first prove (\ref{cor1}).
In the case $k=3$ we  choose $\alpha_3= \alpha_4=\ldots = \alpha_n =0$.
 In this case all $\beta_k$'s are equal to $1/4$.
In the general case $k>3$ we  choose $\alpha_m=-(m-2)/2$, when  $m=3,\ldots,k-1$
and  $\alpha_m=0$,  when  $m=k,\ldots,n$.

To prove (\ref{cor1}) we choose  $\alpha_m=-(m-2)/2$ when  $m=3,\ldots,k-1$,
$a_k=0$, $a_{k+l}= -\frac{l}{2}$, $l=1,...,n-k-1$, $a_n=0$.
\end{proof}

We next  give the proof of the second part of Theorem A: \\

\noindent
{\em Proof of Theorem A, part (ii):}
We will  first  prove that $\beta_3 \leq \frac14$, therefore
  $\beta_3= -\alpha_3^2 +\frac14$,
for suitable  $\alpha_3 \leq 0$. Then, for this $\beta_3$, we will  prove that
$\beta_4 \leq (\alpha_3-\frac12)^2$,   and therefore $\beta_4 = -\alpha_4^2 +
(\alpha_3-\frac12)^2$ for suitable $\xa_4 \leq 0$  and so on.

\noindent
{\bf Step 1.}
Let us first prove  the estimate for  $\beta_3$.  To this end we set
\be\la{22.1}
Q_3[u]:=\frac{\int_{\mathbb{R}^n}  |\nabla u|^2dx - \sum_{i=4}^{n}
 \beta_i \int_{\mathbb{R}^n}
\frac{u^2}{(x_1^2+x_2^2+\ldots+x_i^2)} dx }{\int_{\mathbb{R}^n}
\frac{u^2}{x_1^2+x_2^2+x_3^2} dx }.
\ee
We clearly  have that $\xb_3 \leq \inf_{u \in  C_0^{\infty}(\mathbb{R}^n) } Q_3[u]$.
In the sequel we will show that
\be\la{22.2}
\inf_{u \in  C_0^{\infty}(\mathbb{R}^n) } Q_3[u]
  \leq  \frac14,
\ee
whence,  $\beta_3 \leq \frac14$.

At this point we introduce a family of cutoff functions for later use.
For $j = 3, \ldots, n$ and $k_j > 0$ we set
\[
\phi_j(t) = \left\{ \begin{array}{ll}
0, &  ~~ t < \frac{1}{k_j^2}\\
1+\frac{\ln{k_jt}}{\ln{k_j}}, & ~~ \frac{1}{k_j^2} \leq t < \frac{1}{k_j}\\
1, & ~~ t \geq \frac{1}{k_j},\\
\end{array} \right.
\]
and
\[
h_{k_j}(x) := \phi_j(r_j) \quad \textrm{ where } \quad
r_j:=|\mathbf{X_j}|=  (x_1^2+\ldots+x_j^2)^{\frac{1}{2}}.
\]
Note that
\[
|\nabla h_{k_j}(x)|^2 = \left\{ \begin{array}{ll}
\frac{1}{\ln^2{k_j}}\frac{1}{r_j^2} &  \frac{1}{k_j^2}
 \leq r_j \leq \frac{1}{k_j}\\
0 &  \textrm{otherwise}\\
\end{array} \right..
\]
We also denote by  $ \phi(x)$ a  radially symmetric
 $C_0^{\infty}(\mathbb{R}^n)$ function
such that  $\phi =1$ for $|x|<1/2$
and $\phi =0$ for $|x|>1$.

To prove (\ref{22.2}) we consider the family of functions
\be\la{22.3}
u_{k_3}(x) = |\mathbf{X_3}|^{-\frac12} h_{k_3}(x) \phi(x).
\ee
We will show that as $k_3 \ra \infty$
\be\la{22.4}
\frac{\int_{\mathbb{R}^n}  |\nabla u_{k_3}|^2dx -
 \sum_{i=4}^{n} \beta_i \int_{\mathbb{R}^n}
\frac{u^2_{k_1}}{(x_1^2+x_2^2+\ldots+x_i^2)} dx }{\int_{\mathbb{R}^n}
\frac{u^2_{k_3}}{x_1^2+x_2^2+x_3^2} dx } =
\frac{\int_{\mathbb{R}^n}  |\nabla u_{k_3}|^2dx}{\int_{\mathbb{R}^n}
\frac{u^2_{k_3}}{x_1^2+x_2^2+x_3^2} dx }   + o(1).
\ee
To see this, let us first examine the behavior of the  denominator.
For $k_3$ large we easily compute
\bea\la{22.8a}
\int_{\mathbb{R}^n}
|\mathbf{X_3}|^{-3}h_{k_3}^2\phi^2  dx
& \geq &
C \int_{\frac{1}{k_3} < x_1^2 + x_2^2 +x_3^2< \frac12}
(x_1^2+x_2^2+x_3^2)^{-\frac32} dx_1 dx_2 dx_3 \\
& \geq &  C     \int_{0}^{\pi}    \int_{\frac{1}{k_3}}^{\frac12} r^{-1}\
 \sin \theta dr d \theta \nonumber \\
& \geq &  C  \ln{k_3}. \nonumber
\eea

On the other hand by Lebesgue dominated theorem the terms
$ \sum_{i=4}^{n}
 \beta_i \int_{\mathbb{R}^n_+}
\frac{u^2_{k_3}}{(x_1^2+x_2^2+\ldots+x_i^2)} dx $ are easily seen
to be bounded as $k_3 \ra \infty$.  From this we conclude (\ref{22.4}).

We now estimate the gradient term in (\ref{22.4}).
\be\la{22.4b}
\int_{\mathbb{R}^n}  |\nabla u_{k_3}|^2 dx =
  \frac14  \int_{\mathbb{R}^n}   |\mathbf{X_3}|^{-3}  h_{k_3}^2 \phi^2    dx
+ \int_{\mathbb{R}^n}    |\mathbf{X_3}|^{-1} | \nabla  h_{k_3}|^{2} \phi^2
+ \int_{\mathbb{R}^n}   |\mathbf{X_3}|^{-1}  h_{k_3}^{2} | \nabla \phi|^{2} +  mixed~~ terms.
\ee
The first integral of the right hand side
 behaves exactly as the denominator, that is,
 it goes to infinity
like  $O(\ln{k_3})$. The last
integral  is easily seen to be bounded as  $k_3 \ra \infty$. For the
middle integral  we have
\[
 \int_{\mathbb{R}^n}    |\mathbf{X_3}|^{-1} | \nabla  h_{k_3}|^{2} \phi^2
\leq \frac{C}{\ln^2{k_3}}
 \int_{\frac{1}{k_3^2} \leq (x_1^2+x_2^2+x_3^2)^{1/2} \leq \frac{1}{k_3}}   |\mathbf{X_3}|^{-3} dx_1 dx_2 dx_3
\leq \frac{C}{\ln{k_3}}.
\]
As a consequence of these estimates, we easily get that the mixed terms
 in (\ref{22.4b})
are of the order $o(\ln{k_3})$ as  $k_3 \ra \infty$. Hence, we have that
as  $k_1 \ra \infty$,
\be\la{22.4c}
\int_{\mathbb{R}^n}  |\nabla u_{k_3}|^2 dx =
 \frac14  \int_{\mathbb{R}^n}   |\mathbf{X_3}|^{-3}  h_{k_3}^2 \phi^2    dx
 + o(\ln{k_3}).
\ee
From (\ref{22.4})-(\ref{22.4c})
 we conclude that  as  $k_3 \ra \infty$
\[
  Q_3[u_{k_3}] = \frac14  + o(1),
\]
hence
$\inf_{u \in  C_0^{\infty}(\mathbb{R}^n) } Q_3[u]
  \leq  \frac14$  and consequently  $\xb_3 \leq \frac14$.
Therefore for a suitable nonnegative constant $\alpha_3$ we have that
$\beta_3  =  -\alpha_3^2 + \frac14$. We also set
\be\la{22.5}
\gamma_3 := \alpha_3 +\frac12.
\ee

\noindent
{\bf Step 2.}
We will next show that $\beta_4 \leq  (\alpha_3-\frac12)^2$. To this end,
setting
\be\la{22.6}
Q_4[u] := \frac{\int_{\mathbb{R}^n}  |\nabla u|^2dx - (\frac14-\alpha_3^2)
\int_{\mathbb{R}^n}
\frac{u^2}{x_1^2+x_2^2+x_3^2} dx -
\sum_{i=5}^{n} \beta_i \int_{\mathbb{R}^n}
\frac{u^2}{ |\mathbf{X_i}|^2  } dx }{\int_{\mathbb{R}^n}
\frac{u^2 }{ |\mathbf{X_4}|^2 }dx  },
\ee
will  prove that
\[
\inf_{u \in  C_0^{\infty}(\mathbb{R}^n) } Q_4[u]
  \leq    (\alpha_3-\frac12)^2.
\]
We now consider the family of functions
\bea\la{22.7}
u_{k_3, k_4}(x) & :=  &
|\mathbf{X_3}|^{-\xg_3} |\mathbf{X_4}|^{\alpha_3-\frac12}
 h_{k_3}(x)h_{k_4}(x) \phi(x)  \nonumber \\
&  =: & |\mathbf{X_3}|^{-\xg_3 }v_{k_1, k_2}(x).
\eea
An a easy calculation shows that
\be\la{22.8}
Q_4[u_{k_3,k_4}]=
\frac{\int_{\mathbb{R}^n}|\mathbf{X_3}|^{-2 \xg_3}  |\nabla v_{k_3, k_4} |^2dx -
\sum_{i=5}^{n} \beta_i \int_{\mathbb{R}^n_+}
 |\mathbf{X_3}|^{-2 \xg_3}|\mathbf{X_i}|^{-2}  v_{k_3, k_4}^2  dx }
{\int_{\mathbb{R}^n}
|\mathbf{X_3}|^{- 2\xg_3} |\mathbf{X_4}|^{-2}v_{k_3, k_4}^2  dx  }.
\ee
We  next use the precise form of $v_{k_1,k_2}(x)$. Concerning the
denominator
of $Q_4[u_{k_3,k_4}]$ we have that
\[
\int_{\mathbb{R}^n}
|\mathbf{X_3}|^{- 2\xg_3} |\mathbf{X_4}|^{-2}v_{k_3, k_4}^2  dx
=
\int_{\mathbb{R}^n} (x_1^2+x_2^2+x_3^2)^{-1/2-\alpha_3}
(x_1^2+x_2^2+x_3^2+x_4^2)^{\alpha_1-\frac32} h_{k_3}^2 h_{k_4}^2  \phi^2 dx.
\]
Sending $k_3$ to infinity, using the structure of the cutoff functions and
then introducing polar coordinates we  get
\begin{eqnarray*}
 \int_{\mathbb{R}^n}
|\mathbf{X_3}|^{- 2\xg_3} |\mathbf{X_4}|^{-2}v_{\infty, k_4}^2  dx  =
\int_{\mathbb{R}^n} (x_1^2+x_2^2+x_3^2)^{-1/2-\alpha_3}
(x_1^2+x_2^2+x_3^2+x_4^2)^{\alpha_3-\frac32} h_{k_4}^2  \phi^2 dx,
\end{eqnarray*}
\begin{eqnarray*}
& \geq &
C \int_{\frac{1}{k_4} < x_1^2 + x_2^2 +x_3^2+x_4^2 < \frac12}  (x_1^2+x_2^2+x_3^2)^{-1/2-\alpha_3}
(x_1^2+x_2^2+x_3^2+x_4^2)^{\alpha_3-\frac32}dx_1 dx_2 dx_3 dx_4 \\
& \geq &   C       \int_{\frac{1}{k_4}}^{\frac12} r^{-1}\
  dr \nonumber \\
& \geq &  C  \ln{k_4}. \nonumber
\end{eqnarray*}

The terms  in the numerator   that  are multiplied by the  $\beta_i$'s stay
bounded as $k_3$ or $k_4$ go to infinity.
\bea\la{22.8b}
\int_{\mathbb{R}^n} |\mathbf{X_3}|^{-2 \xg_3}  |\nabla v_{k_3, k_4}|^2 dx &  = &
  \left(\xa_3-\frac12 \right)^2 \int_{\mathbb{R}^n}
  |\mathbf{X_3}|^{-2 \xg_3}  |\mathbf{X_4}|^{2 \xa_3 -3}   h_{k_3}^2  h_{k_4}^2 \phi^2  dx
\nonumber \\
&  & ~+
 \int_{\mathbb{R}^n}  |\mathbf{X_3}|^{-2 \xg_3} |\mathbf{X_4}|^{2 \xa_3 -1}
   |\nabla (h_{k_3} h_{k_4})|^2 \phi^2
 \\
&  & ~+ \int_{\mathbb{R}^n}  |\mathbf{X_3}|^{-2 \xg_3}
  |\mathbf{X_4}|^{2 \xa_3 -1}
  h_{k_3}^2  h_{k_4}^2  |\nabla \phi|^{2}
\nonumber \\
&  & ~+ mixed ~ terms.
    \nonumber
\eea
The first integral in the right hand side above, is the same as the
denominator of $Q_4$, and therefore is finite as $k_3 \ra \infty$ and increases
like $\ln{k_4}$ as $k_4 \ra \infty$, cf (\ref{22.8a}).
 The last integral is bounded, no matter
how big the $k_3$ and $k_4$ are. Concerning the middle term
we have
\bea\la{22.8c}
M[v_{k_3,k_4}] & := &
 \int_{\mathbb{R}^n}  |\mathbf{X_3}|^{-2 \xg_3}
  |\mathbf{X_4}|^{2 \xa_3 -1}
   |\nabla (h_{k_3} h_{k_4})|^2 \phi^2 dx  \nonumber \\
& = &  \int_{\mathbb{R}^n} |\mathbf{X_3}|^{-2 \xg_3}
  |\mathbf{X_4}|^{2 \xa_3 -1}
   |\nabla h_{k_3}|^2  h_{k_4}^2 \phi^2 dx  +
\int_{\mathbb{R}^n} |\mathbf{X_3}|^{-2 \xg_3}
  |\mathbf{X_4}|^{2 \xa_3 -1}
    h_{k_3}^2 |\nabla  h_{k_4}|^2 \phi^2 dx  \nonumber \\
    &+&    ~mixed~~ term  \nonumber \\
& =:&  ~I_1 + I_2 +mixed~~ term.
\eea
Since
\[
 |\mathbf{X_4}|^{2 \xa_3 -1}  h_{k_4}^2 = r_4^{2 \xa_3 -1} \phi_4(r_4)
\leq C_{k_4},~~~~~~~~~~~ 0< r_4 <1,
\]
we easily get
\[
I_1 \leq \frac{C}{(\ln{k_3})^2} \int_{{\frac{1}{k_3^2} < (x_1^2 + x_2^2 +x_3^2)^{1/2} < \frac{1}{k_3}}}
(x_1^2+x_2^2+x_3^2)^{-\alpha_3-\frac{3}{2}} dx_1 dx_2 dx_3,
\]
and therefore, since $\xa_3 \leq 0$,
\be I_1  \leq
\frac{C}{\ln{k_3}}, ~~~~~~~~~~~~~k_3 \ra \infty. \ee Also, since
\[
 |\mathbf{X_3}|^{-2\gamma_3}  h_{k_3}^2 = r_3^{2 \xa_3 -1} \phi_3(r_3)
\leq C_{k_3},~~~~~~~~~~~ 0< r_3 <1,
\]
  we  similarly   get (for any $k_3$)
\bea\la{22.9a}
I_2 &\leq& \frac{C}{(\ln{k_4})^2} \int_{{\frac{1}{k_4^2} < (x_1^2 + x_2^2 +x_3^2+x_4^2)^{1/2} < \frac{1}{k_4}}}
(x_1^2+x_2^2+x_3^2+x_4^2)^{-\frac{1}{2}} dx_1 dx_2 dx_3 dx_4 \\
 &\leq& \frac{C}{\ln{k_4}}, ~~~~~~~~~~k_4 \ra \infty.
\eea

From
(\ref{22.8c})-- (\ref{22.9a}) we have that as $k_4 \ra \infty$,
\[
M[v_{\infty,k_4}] = o(1).
\]
Returning to (\ref{22.8b}) we have that  as $k_4 \ra \infty$,
\be\la{22.10}
\int_{\mathbb{R}^n} |\mathbf{X_3}|^{-2 \xg_3}  |\nabla v_{\infty, k_4}|^2 dx
=   \left(\xa_3-\frac12 \right)^2
\int_{\mathbb{R}^n}
|\mathbf{X_3}|^{- 2 \xg_3} |\mathbf{X_4}|^{-2}v_{\infty, k_4}^2  dx
+  o( \ln{k_4}).
\ee
We then have that as  $k_4 \ra \infty$,
\be\la{22.11}
Q_4[u_{\infty,k_4}]=  \left(\xa_3-\frac12 \right)^2 + o(1),
\ee
consequently,  $\beta_4 \leq \left(\xa_3-\frac12 \right)^2$, and therefore
 $\beta_4 = -\alpha_4^2 +
(\alpha_3-\frac12)^2$ for suitable $\xa_4 \leq 0$. We also set
\[
\gamma_4  = \alpha_4- \alpha_{3} + \frac12.
\]

\noindent
{\bf Step 3.} The general case.
At the ($q-1$)th step,  $3 \leq q \leq n$, we have already established
that
\begin{eqnarray*}
 \beta_3 &=& -\alpha_3^2   +\frac14,  \\
 \beta_m  &=&  - \alpha_m^2 + \left(\alpha_{m-1}-\frac12 \right)^2,
~~~~m=4,5,\ldots,q-1,
\end{eqnarray*}
for suitable nonpositive constants $a_i$. Also, we have defined
\begin{eqnarray*}
 \gamma_3 &=& \alpha_3+ \frac12,  \\
 \gamma_m  &=& \alpha_m- \alpha_{m-1} + \frac12,
~~~~m=4,5,\ldots,q-1.
\end{eqnarray*}
Our goal for the rest of the proof  is to show  that $\beta_q \leq
\left(\alpha_{q-1}-\frac12 \right)^2$. To this end
  we consider  the quotient
\be
Q_q[u] := \frac{\int_{\mathbb{R}^n}  |\nabla u|^2dx
-\sum_{q \neq i=3}^{n} \beta_i \int_{\mathbb{R}^n}
\frac{u^2}{ |\mathbf{X_i}|^2  } dx}
{\int_{\mathbb{R}^n}
\frac{u^2 }{ |\mathbf{X_q}|^2 }dx  }.
\ee
The test function is now given by
\bea
u_{k_3,k_q}(x) & :=  &
|\mathbf{X_3}|^{-\xg_3}   |\mathbf{X_4}|^{-\xg_4}  \ldots
  |\mathbf{X_{q-1}}|^{-\xg_{q-1}}
   |\mathbf{X_q}|^{\alpha_{q-1}-\frac12}
 h_{k_4}(x)
 h_{k_q}(x) \phi(x)  \nonumber \\
&  =: & |\mathbf{X_3}|^{-\xg_3}   |\mathbf{X_4}|^{-\xg_4}
 \ldots   |\mathbf{X_{q-1}}|^{-\xg_{q-1}}
    v_{k_q}(x).
\eea

The proof is analogous to the case $q=4$ and goes along the lines of  \cite{Fil3}.
\finedim

\quad \\
The following corollary is a direct consequence of the above theorem and shows
 that the  constants obtained in Corollary \ref{mainineq} are sharp.

\begin{corollary}\label{co2.3}
For $3 \leq k \leq n$,
\begin{equation} \label{sats1}
\inf_{u \in C^{\infty}_0(\mathbb{R}^n)}
\frac{\int_{\mathbb{R}^n}|\nabla u|^2dx}{\int_{\mathbb{R}^n}
\frac{|u|^2}{|\mathbf{X_k}|^2}} = \left( \frac{k-2}{2} \right)^2,
\end{equation}
\begin{equation} \label{sats2}
\inf_{u \in C^{\infty}_0(\mathbb{R}^n)}
\frac{\int_{\mathbb{R}^n}|\nabla u|^2dx - \left( \frac{k-2}{2} \right)^2
\int_{\mathbb{R}^n}  \frac{|u|^2}{|\mathbf{X_k}|^2} dx-
\frac{1}{4}\int_{\mathbb{R}^n} \frac{|u|^2}{|\mathbf{X_{k+1}}|^2}dx
-\ldots- \frac{1}{4}\int_{\mathbb{R}^n}
\frac{|u|^2}{ |\mathbf{X_m}|^2} dx}{\int_{\mathbb{R}^n}
\frac{|u|^2}{ |\mathbf{X_{m+1}}|^2} dx} = \frac{1}{4}
\end{equation}
for $k\leq m<n$.  And
\begin{equation} \label{sats3}
\inf_{u \in C^{\infty}_0(\mathbb{R}^n)} \frac{
\int_{\mathbb{R}^n}|\nabla u|^2dx - \left( \frac{k-2}{2} \right)^2
\int_{\mathbb{R}^n}\frac{|u|^2}{ |\mathbf{X_k}|^2}dx}{\int_{\mathbb{R}^n}
\frac{|u|^2}{|x|^2}dx} = \left(\frac{n-k}{2}\right)^2.
\end{equation}
\end{corollary}
\begin{proof}
All are consequences of Theorem A.  For  (\ref{sats1})   we   take   $\xa_l=-\frac{l-2}{2}$,
$l=1, \ldots,k-1$.

For  (\ref{sats2}) and  (\ref{sats3})  we take the $a$'s of Corollary \ref{mainineq}.


\end{proof}

\section{Hardy-Sobolev-Maz'ya inequalities}

We first establish the following result that will be used for Theorem B.
\begin{theorem}\la{3.1}
{\bf  (weighted Sobolev inequality)} Let $\xs_2, \xs_3,  \ldots, \xs_n$  be real numbers, with $n
\geq 2$.  We set $c_l :=\xs_2 +\ldots +\xs_l +l-1$, for  $2 \leq l
\leq n$. We assume that
\[
 c_l > 0 ~~~~~~~~~{\rm   whenever} ~~~~~~~~~ \xs_l \neq 0,
\]
for $l=2, \ldots,n$.
 Then, there exists a positive constant $C$ such that
 for any $w   \in C_0^{\infty}(\mathbb{R}^n)$ there holds
\be\la{3.5}
\int_{\mathbb{R}^n}
 |\mathbf{X_2}|^{\xs_2}
 \ldots|\mathbf{X_n}|^{\xs_n} |\nabla w|dx
\geq
C \left(
 \int_{\mathbb{R}^n }   \left(
 |\mathbf{X_2}|^{b} |\mathbf{X_3}|^{\xs_3}
 \ldots|\mathbf{X_n}|^{\xs_n} |w| \right)^{q} dx
 \right)^{\frac{1}{q}},
\ee
where
\[
b=\xs_2-1+\frac{q-1}{q}n \quad \quad   \textrm{ and }\quad \quad 1 < q \leq \frac{n}{n-1}.
\]
\end{theorem}

 {\em Proof:} For
\[
1 < q \leq n/(n-1) \quad \textrm{ and } \quad  b = \xs_2 -1 + \frac{q-1}{q}n,
\]
 we easily obtain the following $L^1$ interpolation inequality
\[
 |||\mathbf{X_2}|^bv||_q \leq c_1 |||\mathbf{X_2}|^{\xs_2}v||_{\frac{n}{n-1}}+c_2|||\mathbf{X_2}|^{\xs_2-1}v||_1.
\]
Using the inequality
\begin{equation}
  \left| \int_{\mathbb{R}^n} \mathrm{div} \mathbf{F} |v|dx \right| \leq \int_{\mathbb{R}^n}| \mathbf{F}||\nabla v|dx, \label{L1}
\end{equation}
with the vector field $\mathbf{F}=|\mathbf{X_2}|^{\xs_2-1}\mathbf{X_2}$ one obtains
\[
 |\xs_2+1| \int_{\mathbb{R}^n}|\mathbf{X_2}|^{\xs_2-1}|v| dx \leq \int_{\mathbb{R}^n} |\mathbf{X_2}|^{\xs_2}|\nabla v|dx.
\]
Here we have to restrict ourselves to $\sigma_2+1 > 0$ to ensure that $|\mathbf{X_2}|^{\xs_2-1}\in L^1_{Loc}(\mathbb{R}^n)$.
Also, by combining this inequality with the standard $L^1$ Sobolev inequality we get
\[
 |||\mathbf{X_2}|^{\xs_2}v||_{\frac{n}{n-1}} \leq |||\mathbf{X_2}|^{\xs_2}|\nabla v|||_1.
\]
Hence we arrive at
\[
  \left(\int_{\mathbb{R}^n}(|\mathbf{X_2}|^{b}|v|)^qdx\right)^{1/q} \leq c \int_{\mathbb{R}^n}|\mathbf{X_2}|^{\xs_2}|\nabla v|dx.
\]
Now let $v=  |\mathbf{X_3}|^{\xs_3}w$ in the above inequality.
This gives
\[
   |||\mathbf{X_2}|^b|\mathbf{X_3}|^{\xs_3}|w|||_q
\leq c \int_{\mathbb{R}^n}|\mathbf{X_2}|^{\xs_2}|\mathbf{X_3}|^{\xs_3}|\nabla w|dx +
 |\xs_3| c \int_{\mathbb{R}^n}|\mathbf{X_2}|^{\xs_2}|\mathbf{X_3}|^{\xs_3-1}|w|dx.
\]
Letting  $\mathbf{F} =
|\mathbf{X_2}|^{\xs_2}|\mathbf{X_3}|^{\xs_3-1}\mathbf{X_3}   $ in (\ref{L1}),
we  get
\be\la{3.11}
 |\xs_2 + \xs_3 + 2|\int_{\mathbb{R}^n}|\mathbf{X_2}|^{\xs_2}|\mathbf{X_3}|^{\xs_3-1}|w|dx
 \leq \int_{\mathbb{R}^n} |\mathbf{X_2}|^{\xs_2}|\mathbf{X_3}|^{\xs_3}|\nabla w|dx.
\ee
Here we have to assume $\xs_2 + \xs_3 + 2 >0$ to guarantee that $|\mathbf{X_2}|^{\xs_2}|\mathbf{X_3}|^{\xs_3-1}\in L^1_{Loc}(\mathbb{R}^n)$.
The two previous estimates give us
\[
  |||\mathbf{X_2}|^b|\mathbf{X_3}|^{\xs_3}|w|||_q
  \leq c  \int_{\mathbb{R}^n} |\mathbf{X_2}|^{\xs_2}|\mathbf{X_3}|^{\xs_3}|\nabla w|dx.
\]
If we would have $\xs_3=0$, we  have our result
immediately and we  do  not  have to check  whether  the constant $\xs_2 +
\xs_3 + 2$ is positive  or not.
\quad
 We may repeat this procedure iteratively. In the $l$-th
step we use the vector field
\[
\mathbf{F} =
|\mathbf{X_2}|^{\xs_2}|\mathbf{X_3}|^{\xs_3}    \ldots
   |\mathbf{X_l}|^{\xs_l-1}\mathbf{X}_l,
\]
in (\ref{L1}) to get :
\[
| c_l|~ |||\mathbf{X_2}|^{\xs_2}|\mathbf{X_3}|^{\xs_3}
  \ldots  |\mathbf{X_l}|^{\xs_l-1} w||_1
 \leq \int_{\mathbb{R}^n} |\mathbf{X_2}|^{\xs_2}|\mathbf{X_3}|^{\xs_3}
   \ldots|\mathbf{X_l}|^{\xs_l} |\nabla w|dx.
\]
As before, we note that we do not need this inequality in the case $\xs_l=0$ and if $\xs_l \neq 0$ we have to assume
$c_l= \xs_2+\ldots + \xs_l +l-1 > 0$ to ensure the integrability of the integrand on the left hand side.
From this it then analogously follows that
\[
 c |||\mathbf{X_2}|^{b}|\mathbf{X_3}|^{\xs_3}
    \ldots  |\mathbf{X_l}|^{\xs_l} w||_q
 \leq \int_{\mathbb{R}^n} |\mathbf{X_2}|^{\xs_2}|\mathbf{X_3}|^{\xs_3}
    \ldots  |\mathbf{X_l}|^{\xs_l} |\nabla w|dx,
\]
which is  (\ref{3.5}). \\

\finedim

For the proof of Theorem C we will use the following variant of Theorem \ref{3.1}.

\begin{theorem}\la{3.2}
({\bf Weighted Sobolev inequality})
 Let $\xs_1, \xs_2,  \ldots, \xs_n$  be real numbers, with $n
\geq 2$.  We set $\overline{c}_l :=\xs_1 +\ldots +\xs_l +l-1$, for
$1 \leq l \leq n$. We assume that
\[
 \overline{c}_l > 0 ~~~~~~~~~{\rm   whenever} ~~~~~~~~~ \xs_l \neq 0,
\]
for $l=1,2, \ldots,n$.
 Then, there exists a positive constant $C$ such that
 for any $w   \in C_0^{\infty}(\mathbb{R}^n)$ there holds
\be\la{3.5b}
\int_{\mathbb{R}^n}|x_1|^{\xs_1}
 |\mathbf{X_2}|^{\xs_2}
 \ldots|\mathbf{X_n}|^{\xs_n} |\nabla w|dx
\geq
C \left(
 \int_{\mathbb{R}^n }   \left(|x_1|^b
 |\mathbf{X_2}|^{\xs_2}
 \ldots|\mathbf{X_n}|^{\xs_n} |w| \right)^{q} dx
 \right)^{\frac{1}{q}},
\ee
where
\[
b=\xs_1-1+\frac{q-1}{q}n \quad \quad   \textrm{ and }\quad \quad 1 < q \leq \frac{n}{n-1}.
\]
\end{theorem}

\noindent {\em Proof:}
Let
\[
1 < q \leq n/(n-1) \quad \textrm{ and } \quad  b = \xs_1 -1 + \frac{q-1}{q}n.
\]
We first consider the case $\xs_1>0$. We will use  the
 following $L^1$ interpolation inequality
\[
 |||x_1|^bv||_q \leq c_1 |||x_1|^{\xs_2}v||_{\frac{n}{n-1}}+c_2|||x_1|^{\xs_2-1}v||_1.
\]
Working similarly as in the proof of Theorem D we end up with
\be\la{E1}
  \left(\int_{\mathbb{R}^n}(|x_1|^{b}|v|)^qdx\right)^{1/q} \leq c \int_{\mathbb{R}^n}|x_1|^{\xs_1}|\nabla v|dx.
\ee
In case $\xs_1=0$, inequality (\ref{E1}) is still valid,  see \cite{Maz}, Section  2.1.6/1.

The rest of the proof goes as in Theorem D. That is, we apply  (\ref{E1})  to $v= |\mathbf{X_2}|^{\xs_3}w$
 to get
\[
   |||x_1|^b|\mathbf{X_2}|^{\xs_2}|w|||_q
\leq c \int_{\mathbb{R}^n}|x_1|^{\xs_1}|\mathbf{X_2}|^{\xs_2}|\nabla w|dx +
 |\xs_2| c \int_{\mathbb{R}^n}|x_1|^{\xs_1}|\mathbf{X_2}|^{\xs_2-1}|w|dx.
\]
Letting  $\mathbf{F} =
|x_1|^{\xs_1}|\mathbf{X_2}|^{\xs_2-1}\mathbf{X_2}   $ in (\ref{L1}),
we  get
\be\la{3.11}
 |\xs_1 + \xs_2 + 1|\int_{\mathbb{R}^n}|x_1|^{\xs_1}|\mathbf{X_2}|^{\xs_2-1}|w|dx
 \leq \int_{\mathbb{R}^n} |x_1|^{\xs_1}|\mathbf{X_2}|^{\xs_2}|\nabla w|dx.
\ee
The condition $\overline{c}_2 = \xs_1 + \xs_2 + 1>0$ guarantees that
 $|x_1|^{\xs_1} |\mathbf{X_2}|^{\xs_2-1} \in L^1_{loc}(\mathbb{R}^n)$ and leads to
\[
   |||x_1|^b|\mathbf{X_2}|^{\xs_2}|w|||_q
\leq c \int_{\mathbb{R}^n}|x_1|^{\xs_1}|\mathbf{X_2}|^{\xs_2}|\nabla w|dx.
\]
We omit further details.

\finedim

We are now ready to give the proof of Theorem B: \\

\noindent
{\em  Proof of Theorem B:}
As a first step we will establish that  for any  $v \in C^{\infty}_0(\mathbb{R}^n )$:
\be\la{3.6}
\int_{\mathbb{R}^n}
 |\mathbf{X_2}|^{2\xs_2-\frac{2QB}{Q+2}}|\mathbf{X_3}|^{\frac{4\xs_3}{Q+2}}
\ldots|\mathbf{X_n}|^{\frac{4\xs_n}{Q+2}}
 |\nabla v|^2 dx
\geq C
 \left( \int_{\mathbb{R}^n}
 |\mathbf{X_2}|^{\frac{2QB}{Q+2}}|\mathbf{X_3}|^{\frac{2Q\xs_3}{Q+2}}
 \ldots|\mathbf{X_n}|^{\frac{2Q\xs_n}{Q+2}}
  |v|^{Q} dx \right)^{\frac{2}{Q}},
\ee
provided that  $c_l := \xs_2 + \ldots + \xs_l + (l-1) >0$, if  $\xs_l \neq 0$,  $2 \leq l\leq n$  where
\[
B=\xs_2-1+\frac{Q-2}{2Q}n ~~~~~~~{\rm and}~~~~~~~2<Q \leq \frac{2n}{n-2}.
\]

To  show (\ref{3.6}) we  apply  Theorem \ref{3.1}  to the function
 $w=|v|^{s}$,  with  $s= \frac{Q+2}{2}$, $sq=Q$ and $b=B$.
Trivial estimates give
\begin{eqnarray*}
 C\left( \int_{\mathbb{R}^n  }|\mathbf{X_2}|^{bq}|\mathbf{X_3}|^{\xs_3q} \ldots
  |\mathbf{X_n}|^{\xs_nq}|v|^{sq}dx \right)^{1/q}
 \leq s \int_{\mathbb{R}^n} |\mathbf{X_2}|^{\xs_3}|\mathbf{X_3}|^{\xs_3}\cdot \ldots \cdot
  |\mathbf{X_n}|^{\xs_n}|v|^{s-1}|\nabla v|dx.
\end{eqnarray*}
 We apply  Cauchy-Schwartz to the
 right hand side and the result  follows.

We will use (\ref{3.6}) with $\xs_2=\frac14((Q-2)n-2Q)$, so that  $2\xs_2-\frac{2QB}{Q+2}=0$.
We notice that  the requirement
\[
c_2 =\xs_2+1=\frac14(Q-2)(n-2)>0,
\]
 is equivalent  to $Q>2$ and therefore is satisfied.

To continue we will use   Lemma  \ref{mainthm}. We recall that
 for $\phi>0$ and $u= \phi v$  with  $v \in C^{\infty}_0(\mathbb{R}^n \setminus S_2)$, we have that
\be\la{3.30}
\int_{\mathbb{R}^n}|\nabla u|^2 dx + \int_{\mathbb{R}^n}\frac{\Delta
\phi}{\phi}|u|^2dx
 =  \int_{\mathbb{R}^n} \phi^2|\nabla v|^2dx.
\ee
We will choose for $\phi$,
\bea
\phi(x) &  =  & |\mathbf{X_3}|^{\frac{2 \xs_3}{Q+2}}  |\mathbf{X_4}|^{\frac{2 \xs_4}{Q+2}}  \ldots
|\mathbf{X_n}|^{\frac{2 \xs_n}{Q+2}} \nonumber \\
 &  =  &
|\mathbf{X_3}|^{-\gamma_3}|\mathbf{X_4}|^{-\gamma_4}\cdot \ldots \cdot
 |\mathbf{X_n}|^{-\gamma_n},
\eea
where,
\begin{eqnarray*}
 \gamma_3 &=& \alpha_3 + \frac12,  \\
 \gamma_m  &=& \alpha_m- \alpha_{m-1} + \frac12,
~~~~m=3,\ldots,n.
\end{eqnarray*}
Therefore
\[
\xs_m = -\frac{Q+2}{2}\gamma_m, \quad~~~~~~~~m=3, \ldots,n.
\]
We now apply  (\ref{3.6}) to obtain that
\be\la{3.18}
 \int_{\mathbb{R}^n} \phi^2|\nabla v|^2dx
\geq C
\left( \int_{\mathbb{R}^n}
|\mathbf{X_2}|^{\frac{Q-2}{2}n-Q}|\phi v|^{Q}dx \right)^{\frac{2}{Q}},
\ee
provided that   for  $3 \leq l \leq n$,
\be\la{3.19}
c_l :=\xs_2 +\ldots +\xs_l +l-1 > 0, ~~~~~
whenever ~~~\xs_l \neq 0.
\ee
Combining (\ref{3.18}) with (\ref{3.30}) we get
\[
\int_{\mathbb{R}^n}|\nabla u|^2 dx + \int_{\mathbb{R}^n}\frac{\Delta
\phi}{\phi}|u|^2dx
 \geq C\left( \int_{\mathbb{R}^n}
|\mathbf{X_2}|^{\frac{Q-2}{2}n-Q}|u|^{Q}dx \right)^{\frac{2}{Q}}.
\]
On the other hand,  by Theorem A(i),
\[
-\frac{\Delta
\phi}{\phi}
=
 \frac{\beta_3}{|\mathbf{X_3}|^2}  + \ldots +
 \frac{\beta_n}{|\mathbf{X_n}|^2},
\]
and the desired inequality follows.

 It remains to check condition (\ref{3.19}). For $l=2$ we have already checked it. For
$3 \leq l \leq n$, after some calculations we find that
\begin{eqnarray*}
c_l & = & \xs_2 +\ldots +\xs_l +l-1  \\
    & = & \frac14(Q-2)(n-2) - \frac{Q+2}{2}(\gamma_3   +  \ldots + \gamma_l) + l-1 \\
    & = & \frac{Q+2}{2} \left( -\alpha_l + \frac{(Q-2)(n-l)}{2(Q+2)} \right)  .
\end{eqnarray*}
Recalling that   $\alpha_l \leq 0$, we conclude that   if $l \leq n-1$ then $c_l >0$, whereas if $l=n$, then
$c_n >0$ if and only if $\alpha_n < 0$.
  This proves  (\ref{1.20}) for   $u \in C^{\infty}_0(\mathbb{R}^n \setminus S_2)$  and by a density argument
the result holds for any  $u \in C^{\infty}_0(\mathbb{R}^n)$

In the rest of the proof we will show that (\ref{1.20}) fails in case $\xa_n=0$.
To this end we will establish that
\be\la{3.25}
\inf_{u \in C^{\infty}_0(\mathbb{R}^n)} \frac{
\int_{\mathbb{R}^n}|\nabla u|^2dx - \xb_3
\int_{\mathbb{R}^n}\frac{|u|^2}{ |\mathbf{X_{3}}|^2}dx - \ldots - \xb_n
\int_{\mathbb{R}^n}\frac{|u|^2}{ |\mathbf{X_{n}}|^2}dx}
{\left(\int_{\mathbb{R}^n}|\mathbf{X_{2}}|^{\frac{Q-2}{2}n-Q}|u|^{Q}dx
  \right)^{\frac{2}{Q}}}=0,
\ee
where  $\xb_n =  \left(\alpha_{n-1}-\frac12 \right)^2$. Let
\[
u(x)=
  |\mathbf{X_3}|^{-\xg_3}
 \ldots   |\mathbf{X_{n-1}}|^{-\xg_{n-1}}
    v(x).
\]
A straightforward calculation, quite similar to the one leading to (\ref{22.8}), shows
that  the infimum in (\ref{3.25}) is the same as the following  infimum
\be\la{3.60}
\inf_{v \in C^{\infty}_0(\mathbb{R}^n)}
\frac{\int_{\mathbb{R}^n}
\prod_{j=3}^{n-1}  |\mathbf{X_{j}}|^{-2 \xg_{j}}
 |\nabla v |^2dx -
 \beta_n \int_{\mathbb{R}^n}
\prod_{j=3}^{n-1}  |\mathbf{X_{j}}|^{-2 \xg_{j}}
|\mathbf{X_n}|^{-2}  v ^2  dx}
{\left( \int_{\mathbb{R}^n}
\left(|\mathbf{X_2}|^{\frac{Q-2}{2Q}n-1} \prod_{j=3}^{n-1}|\mathbf{X_{j}}|^{- \xg_{j}} \right)^{Q}
  |v|^{Q}  dx \right)^{\frac{2}{Q}}}.
\ee
We  now choose the following test functions
\be\la{tf}
v_{k_3, \xe}=  |\mathbf{X_{n}}|^{- \xg_{n} +\xe} h_{k_3}(x) \phi(x),~~~~~~~~~~~~~\xe>0,
\ee
where $ h_{k_3}(x)$ and $\phi(x)$ are the same test functions as in  the first step
of the proof of Theorem  A(ii). For this choice, after straightforward
calculations,  quite similar to the ones used in the proof of Theorem A(ii),
we obtain the following estimate for the numerator $N$ in (\ref{3.60}).
\begin{eqnarray*}
N[v_{\infty, \xe}] &  =  & \left( \left(\alpha_{n-1}-\frac12  +\xe \right)^2-
\left(\alpha_{n-1}-\frac12  \right)^2 \right)
\int_{\mathbb{R}^n}
\prod_{j=3}^{n-1}  |\mathbf{X_{j}}|^{-2 \xg_{j}}
 |\mathbf{X_{n}}|^{- 2 \xg_{n}+2 +\xe} \phi^2(x) dx + O_{\xe}(1),  \\
& =  &
 C  \xe  \int_{\mathbb{R}^n}
 r^{-1+ 2 \xe}  \sin{\theta_2}\prod_{j=3}^{n-1} ( \sin{\theta_j} )^{1-2 \xa_j}  \phi^2(r)
 d \theta_1 \ldots d \theta_{n-1}  d r + O_{\xe}(1)    \\
& =  & C \xe  \int_{0}^{1} r^{-1+\xe} dr + O_{\xe}(1).
\end{eqnarray*}
In the above calculations we have taken the limit $k_3 \ra \infty$
and we have used polar coordinates in $(x_1,\ldots, x_n) \ra
(\theta_1,\ldots, \theta_{n-1},r)$. We then conclude that
\be\la{3.81} N[v_{\infty, \xe}]  < C, ~~~~~~~~~~~{\rm as}
~~~~~~~~\xe \ra 0. \ee Similar calculations for the denominator
$D$  in (\ref{3.60}) reveal that
\begin{eqnarray*}
D[v_{\infty, \xe}] &  =  & C
\left(\int_{\mathbb{R}^n}  r^{-1+\xe Q}
\prod_{j=2}^{n-1}(\sin{\theta_j})^{j-n-1+Q(\frac{n-j}{2}-\alpha_j)}
\phi^{Q}
 d\theta_1 \ldots d\theta_{n-1} dr   \right)^{\frac{2}{Q}}  \\
 &  \geq   & C \left( \int_{0}^{\frac12}
 r^{-1+\xe Q} dr \right)^{\frac{2}{Q}} \\
& = & C \xe^{- \frac{2}{Q}}.
\end{eqnarray*}
We then  have that
\[
\frac{N[v_{\infty, \xe}]}{D[v_{\infty, \xe}] } \ra 0 ~~~~~~{\rm as}~~~~
\xe \ra 0,
\]
and therefore the  infimum in (\ref{3.60})  or (\ref{3.25}) is equal to zero.
This completes the proof of the Theorem.

\finedim

Here is a consequence of Theorem  B.

\begin{corollary}
Let  $3 \leq k<n$  and  $2<Q \leq \frac{2n}{n-2}$.  Then, for any
$\beta_n < \frac{1}{4}$, there exists a positive constant $C$ such
that for  all  $u \in C^{\infty}_0(\mathbb{R}^n)$ there holds
\begin{eqnarray*}
  \int_{\mathbb{R}^n}|\nabla u|^2dx &\geq& \int_{\mathbb{R}^n}
 \left( \left(\frac{k-2}{2}\right)^2 \frac{1}{ |\mathbf{X_k}|^2}
 + \frac{1}{4}\frac{1}{ |\mathbf{X_{k-1}}|^2}
 +\ldots \right. + \left. \frac{1}{4}\frac{1}{ |\mathbf{X_{n-1}}|^2}+
 \frac{\beta_n}{|\mathbf{X_n}|^2}
 \right)|u|^2dx\\
  &+& C\left(\int_{\mathbb{R}^n} |\mathbf{X_2}|^{\frac{Q-2}{2}n-Q}|u|^{Q}dx
  \right)^{\frac{2}{Q}}.
\end{eqnarray*}
If $\beta_n = \frac{1}{4}$ the previous inequality fails.

In case  $k=n$ we have that for  any  $\beta_n < \frac{(n-2)^2}{4}$,
 there exists a positive constant $C$
such that for  all  $u \in C^{\infty}_0(\mathbb{R}^n)$ there holds
\[
 \int_{\mathbb{R}^n}|\nabla u|^2dx \geq
 \xb_n\int_{\mathbb{R}^n}\frac{u^2}{|x|^2}dx +
 C\left(\int_{\mathbb{R}^n}|\mathbf{X_2}|^{\frac{Q-2}{2}n-Q}|u|^{Q}dx
  \right)^{\frac{2}{Q}}.
\]
The above inequality fails for $\beta_n = \frac{(n-2)^2}{4}$.
\end{corollary}
\begin{proof} In Theorem B we make the following choices:
In the case $k=3$ we  choose $\alpha_3= \alpha_4=\ldots = \alpha_{n-1} =0$.
 In this case  $\beta_k=1/4$,    $k=1,\ldots,n-1$. The condition $\xa_n<0$ is
equivalent to $\beta_n < \frac{1}{4}$.

In the case  $3<k \leq n-1$ we  choose $\alpha_m=-m/2$, when  $m=1,2,\ldots,k-1$
and  $\alpha_m=0$,  when  $m=k,\ldots,n-1$.
Finally, in case  $k=n$,  we  choose $\alpha_m=-(m-2)/2$, for  $m=3,4,\ldots,n-1$.

\end{proof}

We finally give the proof of Theorem C:  \\
{\em Proof of Theorem C:}
We first prove that  the following   inequality holds for any    $v \in C^{\infty}_0(\mathbb{R}^n )$:
\be\la{3.c6}
\int_{\mathbb{R}^n}
 |x_1|^{2\xs_1-\frac{2QB}{Q+2}}|\mathbf{X_2}|^{\frac{4\xs_2}{Q+2}}
\ldots|\mathbf{X_n}|^{\frac{4\xs_n}{Q+2}}
 |\nabla v|^2 dx
\geq C
 \left( \int_{\mathbb{R}^n}
 |x_1|^{\frac{2QB}{Q+2}}|\mathbf{X_2}|^{\frac{2Q\xs_2}{Q+2}}
 \ldots|\mathbf{X_n}|^{\frac{2Q\xs_n}{Q+2}}
  |v|^{Q} dx \right)^{\frac{2}{Q}},
\ee
provided that  $\overline{c}_l := \xs_1 + \ldots + \xs_l + (l-1) >0$, if  $\xs_l \neq 0$,  $1 \leq l\leq n$  where
\[
B=\xs_1-1+\frac{Q-2}{2Q}n ~~~~~~~{\rm and}~~~~~~~ \frac{2(n-1)}{n-2}<Q \leq \frac{2n}{n-2}.
\]

To  show (\ref{3.c6}) we  apply Theorem  \ref{3.2} to the function
 $w=|v|^{s}$,  with  $s= \frac{Q+2}{2}$, $sq=Q$ and $b=B$, and then use  Cauchy-Schwartz inequality.

We will use (\ref{3.c6}) with $\xs_1=\frac14((Q-2)n-2Q)$   and   $\xs_2=0$. In this case   $2\xs_1-\frac{2QB}{Q+2}=0$.
The choice of  $\phi$  stays  the same as  in the proof of Theorem B.  Eventually, we arrive at
\[
\int_{\mathbb{R}^n}|\nabla u|^2 dx - \int_{\mathbb{R}^n} \left(  \frac{\beta_3}{|\mathbf{X_3}|^2}  + \ldots +
 \frac{\beta_n}{|\mathbf{X_n}|^2} \right) |u|^2dx
 \geq C\left( \int_{\mathbb{R}^n}
|x_1|^{\frac{Q-2}{2}n-Q}|u|^{Q}dx \right)^{\frac{2}{Q}},
\]
provided that  the  $\overline{c}_l$'s satisfy our assumptions of
Theorem \ref{3.2}.  However it turns out that
\[
\overline{c}_l =  \frac{Q+2}{2} \left( -\alpha_l + \frac{(Q-2)(n-l)}{2(Q+2)} \right),
~~~~~~~~1\leq l \leq n,
\]
and our assumptions are satisfied in case $\alpha_n<0$.

In remains to prove  that (\ref{1.20a}) fails in case $\xa_n=0$.
To this end we will establish that
\be\la{3.25a}
\inf_{u \in C^{\infty}_0(\mathbb{R}^n)} \frac{
\int_{\mathbb{R}^n}|\nabla u|^2dx - \xb_3
\int_{\mathbb{R}^n}\frac{|u|^2}{ |\mathbf{X_{3}}|^2}dx - \ldots - \xb_n
\int_{\mathbb{R}^n}\frac{|u|^2}{ |\mathbf{X_{n}}|^2}dx}
{\left(\int_{\mathbb{R}^n}|x_1|^{\frac{Q-2}{2}n-Q}|u|^{Q}dx
  \right)^{\frac{2}{Q}}}=0,
\ee
where  $\xb_n =  \left(\alpha_{n-1}-\frac12 \right)^2$.
The test functions used  in the proof of Theorem B  can  also be used here since they belong in the proper function
 space.  The result follows by observing that the weight here is stronger than in Theorem B.

\finedim

An easy consequence of the above Theorem is the following:

\begin{corollary}
Let  $3 \leq k<n$  and  $\frac{2(n-1)}{n-2} <  Q  \leq \frac{2n}{n-2}$.  Then, for any
$\beta_n < \frac{1}{4}$, there exists a positive constant $C$ such
that for  all  $u \in C^{\infty}_0(\mathbb{R}^n)$ there holds
\begin{eqnarray*}
  \int_{\mathbb{R}^n}|\nabla u|^2dx &\geq& \int_{\mathbb{R}^n}
 \left( \left(\frac{k-2}{2}\right)^2 \frac{1}{ |\mathbf{X_k}|^2}
 + \frac{1}{4}\frac{1}{ |\mathbf{X_{k-1}}|^2}
 +\ldots \right. + \left. \frac{1}{4}\frac{1}{ |\mathbf{X_{n-1}}|^2}+
 \frac{\beta_n}{|\mathbf{X_n}|^2}
 \right)|u|^2dx\\
  &+& C\left(\int_{\mathbb{R}^n} |x_1|^{\frac{Q-2}{2}n-Q}|u|^{Q}dx
  \right)^{\frac{2}{Q}}.
\end{eqnarray*}
If $\beta_n = \frac{1}{4}$ the previous inequality fails.

In case  $k=n$ we have that for  any  $\beta_n < \frac{(n-2)^2}{4}$,
 there exists a positive constant $C$
such that for  all  $u \in C^{\infty}_0(\mathbb{R}^n)$ there holds
\[
 \int_{\mathbb{R}^n}|\nabla u|^2dx \geq
 \xb_n\int_{\mathbb{R}^n}\frac{u^2}{|x|^2}dx +
 C\left(\int_{\mathbb{R}^n}|x_1|^{\frac{Q-2}{2}n-Q}|u|^{Q}dx
  \right)^{\frac{2}{Q}}.
\]
The above inequality fails for $\beta_n = \frac{(n-2)^2}{4}$.
\end{corollary}

\end{document}